\documentclass[12pt,a4paper]{amsart}
\usepackage{amssymb,amsfonts,amsthm,amsmath,graphicx,xcolor,csquotes}
\usepackage[scr]{rsfso}
\usepackage{chngcntr}
\usepackage{enumitem}

\theoremstyle{plain}
\newtheorem{theorem}{Theorem}
\newtheorem{lemma}{Lemma}

\newtheorem{corollary}{Corollary}

\numberwithin{theorem}{section}
\numberwithin{equation}{section}
\numberwithin{statement}{section}
\numberwithin{lemma}{section}
\numberwithin{definition}{section}
\numberwithin{conjecture}{section}
\numberwithin{corollary}{section}
\begin{document}
\dedicatory{Dedicated to Professor Rodney J Baxter on his $83$rd birthday.}
\title[continued fraction partition identities]{Continued Fractions for partition generating functions}
\author{Geoffrey B Campbell}
\address{Mathematical Sciences Institute,
         The Australian National University,
         Canberra, ACT, 0200, Australia}

\email{Geoffrey.Campbell@anu.edu.au}

\thanks{Thanks are due to Professor Dr Henk Koppelaar, whose discussions and suggestions have been very helpful for the book for which this paper is essentially a chapter.}

\keywords{Continued fractions and generalizations. Exact enumeration problems, generating functions. Partitions of integers. Elementary theory of partitions. Combinatorial identities, bijective combinatorics. Lattice points in specified regions.}
\subjclass[2010]{Primary: 11J70; Secondary: 05A15, 05E40, 11Y11, 11P21}

\begin{abstract}
We derive continued fractions for partition generating functions, utilizing both Euler's techniques and Ramanujan's techniques. Although our results are for integer partitions there is scope to extend this work to vector partitions, including for binary and n-ary partitions.
\end{abstract}

\maketitle

\section{Euler's Continued Fraction}

Almost 290 years ago in 1737, Leonhard Euler wrote \textit{De fractionibus continuis dissertatio}, which gave mathematics a first ever comprehensive account of the properties of continued fractions, and included the first proof that the number $e$ is irrational. (See Sandifer \cite{cS2006})  Later, but still 275 years ago in 1748, Euler, in his \textit{Introductio in analysin infinitorum Vol. I, Chapter 18} \cite{lE1748}, proved:
\begin{enumerate}
  \item the equivalence of his continued fraction to a generalized infinite series,
  \item every rational number can be written as a finite continued fraction, and
  \item the continued fraction of an irrational number is infinite.
\end{enumerate}

Euler’s continued fraction is the very nice identity, whose first few cases are:
\begin{eqnarray*}
  a_0 + a_0 a_1  &=& a_0/(1 - a_1/(1 + a_1)) \\
                 &=& \cfrac{a_0}{1- \cfrac{a_1}{1+a_1}} \; ;  \\
  a_0 + a_0 a_1 + a_0 a_1 a_2  &=& a_0/(1 - a_1/(1 + a_1 - a_2/(1 + a_2))) \\
                 &=& \cfrac{a_0}{1- \cfrac{a_1}{1+a_1-\cfrac{a_2}{1+a_2}}} \; ;  \\
  a_0 + a_0 a_1 + a_0 a_1 a_2 + a_0 a_1 a_2 a_3  &=& a_0/(1 - a_1/(1 + a_1 - a_2/(1 + a_2 - a_3/(1 + a_3)))) \\
                 &=& \cfrac{a_0}{1- \cfrac{a_1}{1+a_1-\cfrac{a_2}{1+a_2-\cfrac{a_3}{1+a_3 }}}} \; .
\end{eqnarray*}

Hence, we can state Euler's Continued Fraction in the following

\begin{theorem}  \label{cf.01}
If $a_0$, $a_1$, $a_3$, ... $a_n$ are defined functions such that no denominator is zero in the following equations then
  \begin{equation} \label{cf.02}
  \sum_{k=0}^{n} \prod_{j=0}^{k} a_j = a_0 + a_0 a_1 + a_0 a_1 a_2 + ... + a_0 a_1...a_n
\end{equation}
\begin{equation*}
= a_0/(1 - a_1/(1 + a_1 - a_2/(1 + a_2 - a_3/(1 + ... \, a_{n-1}/(1 + a_{n-1} - a_n/(1 + a_n))))).
\end{equation*}
\begin{equation*}
   = \cfrac{a_0}{1- \cfrac{ a_1}{1 + a_1 -\cfrac{a_2}{1 + a_2 -\cfrac{a_3}{1 + a_3 -\cfrac{\ddots}{\ddots \cfrac{a_{n-1}}{1+ a_{n-1} -\cfrac{a_n}{1 + a_n}}}}}}} \; .
\end{equation*}
\end{theorem}

Obviously, this lends itself to many of the elementary series that arise in school and university analysis. However, we shall put this to good use in applying it to partition generating functions. The fact of this theorem involving a finite sum allows us to incrementally extend the number of terms until we can infer the infinite versions of the theorem.

Example 1: The exponential function is
  \begin{equation} \label{cf.03}
  \exp(z) = 1 + \frac{z}{1!} + \frac{z^2}{2!} + \frac{z^3}{3!} + ... = 1 + \left(\frac{z}{1}\right) + \left(\frac{z}{1}\right)\left(\frac{z}{2}\right) + \left(\frac{z}{1}\right)\left(\frac{z}{2}\right)\left(\frac{z}{3}\right) + ...
\end{equation}
\begin{equation*}
  = 1/\left(1 - z/\left(1 + z - \left(\frac{z}{2}\right)/\left(1 + \left(\frac{z}{2}\right) - \left(\frac{z}{3}\right)/\left(1 + \left(\frac{z}{3}\right) - \left(\frac{z}{4}\right)/\left(1 + \left(\frac{z}{4}\right) - ... \right)\right)\right)\right)\right).
\end{equation*}
Applying an “equivalence transformation” that consists of clearing the fractions, this example is simplified to
\begin{equation*}
  \exp(z) = 1/(1 - z/(1 + z - z/(2 + z - 2z/(3 + z - 3z/(4 + z - \ldots ))))),
\end{equation*}
or the equivalent statement
\begin{equation*}
 \exp(z) = \cfrac{1}{1-\cfrac{z}{1+z-\cfrac{z}{2+z- \cfrac{2z}{3+z- \cfrac{3z}{4+z- \dots}}}}}
\end{equation*}
and we know this continued fraction converges uniformly on every bounded domain in the complex plane because it is equivalent to the power series for $\exp(z)$.

Example 2: There is the well-known logarithmic function series
  \begin{equation} \label{cf.04}
\log\left(\frac{1+z}{1-z}\right)= 2z (\frac{1}{1} + \frac{z^2}{3} + \frac{z^4}{5} + ... )
\end{equation}
\begin{equation*}
= 2z (1 + (\frac{z^2}{3}) + (\frac{z^2}{3})(\frac{3z^2}{5}) + (\frac{z^2}{3})(\frac{3z^2}{5})(\frac{5z^2}{7}) + ... ).
\end{equation*}
Applying Euler's continued fraction formula to this expression shows that:
\begin{equation*}
\log\left(\frac{1+z}{1-z}\right)
\end{equation*}
\begin{equation*}
= 2z/(1 - (\frac{z^2}{3})/(1 + (\frac{z^2}{3}) - (\frac{3z^2}{5})/(1 + (\frac{3z^2}{5}) - (\frac{5z^2}{7})/(1 + (\frac{5z^2}{7}) - (\frac{7z^2}{9})/(1 + (\frac{7z^2}{9}) - ... ))))).
\end{equation*}
Applying the “equivalence transformation” this example is simplified to
\begin{equation*}
\log\left(\frac{1+z}{1-z}\right) = 2z/(1-z^2/(z^2 + 3 -(3z)^2/(3z^2 + 5 -(5z)^2/( 5z^2 + 7 -(7z)^2/( 7z^2 + 9 - ... )))))
\end{equation*}
\begin{equation*}
  = \cfrac{2z}{1-\cfrac{z^2}{z^2 +3-\cfrac{(3z)^2}{3z^2 +5- \cfrac{(5z)^2}{5z^2 +7- \cfrac{(7z)^2}{7z^2 + 9- \dots}}}}}
\end{equation*}

Example 3: A continued fraction for $\pi$.
We can use the previous example involving the principal branch of the natural logarithm function to construct a continued fraction representation of $\pi$. First we note that
\begin{equation*}
(i+1)/(i-1) = i, \quad so \quad then \quad \log((i+1)/(i-1)) = i \pi/2.
\end{equation*}
Setting $z = i$ in the previous result, and remembering that $i^2 = -1$, we obtain immediately
\begin{equation*}
  \pi = \cfrac{4}{1+ \cfrac{1^2}{2 +\cfrac{3^2}{2 + \cfrac{5^2}{2 + \cfrac{7^2}{2 + \dots}}}}}
\end{equation*}

\section{Euler's continued fraction applied to partitions}

In this section we will technically do no more than apply the previous section. However, the theory of partitions is full of generating functions that are emenable to the Euler continued fraction. In a subsequent section we will examine Ramanujan type continued fractions, but firstly we will gather some "low hanging fruit" from some elementary series-product identities.

We begin with the well-known telescoping identities:

If $a_1$, $a_2$, $a_3$, ... , $a_n$, are functions chosen for nonzero denominators, then
\begin{equation} \label{cf.05}
  1 + \frac{a_1}{1-a_1} + \frac{a_2}{(1-a_1)(1-a_2)} + ... + \frac{a_n}{(1-a_1)(1-a_2)...(1-a_n)}
\end{equation}
\begin{equation*}
 = \frac{1}{(1-a_1)(1-a_2)(1-a_3)...(1-a_n)};
\end{equation*}
and
\begin{equation} \label{cf.06}
  1 + a_1+ a_2(1+a_1) + a_3(1+a_1)(1+a_2) + ... + a_n(1+a_1)(1+a_2)...(1+a_{n-1})
\end{equation}
\begin{equation*}
 = (1+a_1)(1+a_2)(1+a_3)...(1+a_n).
\end{equation*}

The series in (\ref{cf.05}) and (\ref{cf.06}) are already close to being in the required form to apply the Euler continued fraction since
\begin{equation} \label{cf.07}
  1 + \frac{a_1}{1-a_1} + \frac{a_2}{(1-a_1)(1-a_2)} + ... + \frac{a_n}{(1-a_1)(1-a_2)...(1-a_n)}
\end{equation}
\begin{equation*}
 = 1 + \frac{a_1}{1-a_1} + \frac{a_1}{1-a_1}\frac{a_2(1-a_1)}{a_1(1-a_2)} + ... + \frac{a_1}{1-a_1}\frac{a_2(1-a_1)}{a_1(1-a_2)}...\frac{a_n(1-a_{n-1})}{a_{n-1}(1-a_n)};
\end{equation*}
and
\begin{equation} \label{cf.08}
  1 + a_1+ a_2(1+a_1) + a_3(1+a_1)(1+a_2) + ... + a_n(1+a_1)(1+a_2)...(1+a_{n-1})
\end{equation}
\begin{equation*}
 = 1 + a_1+ a_1 \frac{a_2(1+a_1)}{a_1} + a_1 \frac{a_2(1+a_1)}{a_1}\frac{a_3(1+a_2)}{a_2} + ... + a_1 \frac{a_2(1+a_1)}{a_1}\frac{a_3(1+a_2)}{a_2}...\frac{a_n(1+a_{n-1})}{a_{n-1}}.
\end{equation*}

Hence combining (\ref{cf.05}) with (\ref{cf.07}) and then (\ref{cf.06}) with (\ref{cf.08}) respectively, we obtain

  \begin{equation} \label{cf.09}
  \frac{1}{(1-a_1)(1-a_2)(1-a_3)...(1-a_n)}
\end{equation}
\begin{equation*}
   = \cfrac{1}{1- \cfrac{\frac{a_1}{1-a_1}}{1 + \frac{a_1}{1-a_1} -\cfrac{\frac{a_2(1-a_1)}{a_1(1-a_2)}}{1 + \frac{a_2(1-a_1)}{a_1(1-a_2)} -\cfrac{\frac{a_3(1-a_2)}{a_2(1-a_3)}}{1 + \frac{a_3(1-a_2)}{a_2(1-a_3)} -\cfrac{\ddots}{\ddots \cfrac{\frac{a_{n-1}(1-a_{n-2})}{a_{n-2}(1-a_{n-1})}}{1+ \frac{a_{n-1}(1-a_{n-2})}{a_{n-2}(1-a_{n-1})} -\cfrac{\frac{a_n(1-a_{n-1})}{a_{n-1}(1-a_n)}}{1 + \frac{a_n(1-a_{n-1})}{a_{n-1}(1-a_n)}}}}}}}} \; ;
\end{equation*}
and
  \begin{equation} \label{cf.10}
  (1+a_1)(1+a_2)(1+a_3)...(1+a_n)
\end{equation}
\begin{equation*}
   = \cfrac{1}{1- \cfrac{a_1}{1 + a_1 -\cfrac{\frac{a_2(1+a_1)}{a_1}}{1 + \frac{a_2(1+a_1)}{a_1} -\cfrac{\frac{a_3(1+a_2)}{a_2}}{1 + \frac{a_3(1+a_2)}{a_2} -\cfrac{\ddots}{\ddots \cfrac{\frac{a_{n-1}(1+a_{n-2})}{a_{n-2}}}{1+ \frac{a_{n-1}(1+a_{n-2})}{a_{n-2}}
   -\cfrac{\frac{a_n(1+a_{n-1})}{a_{n-1}}}{1 + \frac{a_n(1+a_{n-1})}{a_{n-1}}}}}}}}} \; .
\end{equation*}

After applying the “equivalence transformation” to both of (\ref{cf.09}) and then (\ref{cf.10}) to eliminate denominator terms, each continued fraction is simplified giving us the following two theorems.

\begin{theorem} \label{cf.11}
If $a_1$, $a_2$, $a_3$, ... , $a_n$, are functions chosen for nonzero denominators, then
    \begin{equation} \label{cf.12}
  \frac{1}{(1-a_1)(1-a_2)(1-a_3)...(1-a_n)}
\end{equation}
\begin{equation*}
   = \cfrac{1}{1- \cfrac{a_1}{1 -\cfrac{a_2}{a_1 +a_2 - 2 a_1 a_2 -\cfrac{a_1 a_3}{a_2 + a_3 - 2 a_2 a_3 -\cfrac{\ddots}{\ddots
   \cfrac{a_{n-2} a_n}{a_{n-1} + a_n - 2a_{n-1}a_n }}}}}} \; .
\end{equation*}
\end{theorem}

At first glance we can see this theorem as being applicable to generating functions for unrestricted partitions of various kinds. Similarly the next theorem applies for partitions of various sorts into distinct parts.

\begin{theorem} \label{cf.13}
If $a_1$, $a_2$, $a_3$, ... , $a_n$, are functions chosen for nonzero denominators, then
    \begin{equation} \label{cf.14}
  (1+a_1)(1+a_2)(1+a_3)...(1+a_n)
\end{equation}
\begin{equation*}
   = \cfrac{1}{1- \cfrac{a_1}{1 + a_1 -\cfrac{(1+a_1)a_2}{a_1 + a_2 + a_1 a_2 -\cfrac{(1+a_2)a_3}{a_2 + a_3 + a_2 a_3 -\cfrac{\ddots}{\ddots \cfrac{(1+a_{n-1})a_n}{a_{n-1}+a_n+a_{n-1}a_n}}}}}} \; .
\end{equation*}
\end{theorem}

There are many examples we could choose for substitution into theorems \ref{cf.13} and \ref{cf.13}. So, let's start with the generating functions for unrestricted partitions, and for distinct partitions as follows.

\begin{corollary} \label{cf.15}
If $p_n(k)$, is the number of unrestricted partitions of $k$ into integers no greater than $n$, then
    \begin{equation} \label{cf.16}
  \frac{1}{(1-q^1)(1-q^2)(1-q^3)...(1-q^n)} = \sum_{k=0}^{\infty} p_n(k) q^k
\end{equation}
\begin{equation*}
   = \cfrac{1}{1- \cfrac{q^1}{1 -\cfrac{q^2}{q^1 +q^2 - 2 q^1 q^2 -\cfrac{q^1 q^3}{q^2 + q^3 - 2 q^2 q^3 -\cfrac{\ddots}{\ddots
   \cfrac{q^{n-2} q^n}{q^{n-1} + q^n - 2q^{n-1}q^n }}}}}} \; .
\end{equation*}
\end{corollary}

\begin{corollary} \label{cf.17}
If $p_n(\mathfrak{D},k)$, is the number of distinct partitions of $k$ into integers no greater than $n$, then
    \begin{equation} \label{cf.18}
  (1+q^1)(1+q^2)(1+q^3)...(1+q^n) = \sum_{k=0}^{\infty} p_n(\mathfrak{D},k) q^k
\end{equation}
\begin{equation*}
   = \cfrac{1}{1- \cfrac{q^1}{1 + q^1 -\cfrac{(1+q^1)q^2}{q^1 + q^2 + q^1 q^2 -\cfrac{(1+q^2)q^3}{q^2 + q^3 + q^2 q^3 -\cfrac{\ddots}{\ddots \cfrac{(1+q^{n-1})q^n}{q^{n-1}+q^n+q^{n-1}q^n}}}}}} \; .
\end{equation*}
\end{corollary}

Next we choose the odd integer powers substituted into the two theorems.

\begin{corollary} \label{cf.19}
If $p_n(\mathfrak{O},k)$, is the number of unrestricted partitions of $k$ into odd integers no greater than $2n-1$, then
    \begin{equation} \label{cf.20}
  \frac{1}{(1-q^1)(1-q^3)(1-q^5)...(1-q^{2n-1})} = \sum_{k=0}^{\infty} p_n(\mathfrak{O},k) q^k
\end{equation}
\begin{equation*}
   = \cfrac{1}{1- \cfrac{q^1}{1 -\cfrac{q^3}{q^1 +q^3 - 2 q^1 q^3 -\cfrac{q^1 q^5}{q^3 + q^5 - 2 q^3 q^5 -\cfrac{\ddots}{\ddots
   \cfrac{q^{n-2} q^n}{q^{2n-3} + q^{2n-1} - 2q^{2n-3}q^{2n-1}}}}}}} \; .
\end{equation*}
\end{corollary}

\begin{corollary} \label{cf.21}
If $p_n(\mathfrak{DO},k)$, is the number of distinct partitions of $k$ into odd integers no greater than $2n-1$, then
    \begin{equation} \label{cf.22}
  (1+q^1)(1+q^3)(1+q^5)...(1+q^{2n-1}) = \sum_{k=0}^{\infty} p_n(\mathfrak{DO},k) q^k
\end{equation}
\begin{equation*}
   = \cfrac{1}{1- \cfrac{q^1}{1 + q^1 -\cfrac{(1+q^1)q^3}{q^1 + q^3 + q^1 q^3 -\cfrac{(1+q^3)q^5}{q^3 + q^5 + q^3 q^5 -\cfrac{\ddots}{\ddots \cfrac{(1+q^{2n-3})q^{2n-1}}{q^{2n-3}+q^{2n-1}+q^{2n-3}q^{2n-1}}}}}}} \; .
\end{equation*}
\end{corollary}

It is a well-known result due to Euler that $p_\infty(\mathfrak{DO},k) = p_\infty(\mathfrak{O},k)$. Explicitly, as $n \rightarrow \infty$ equations (\ref{cf.22}) and (\ref{cf.20}) are equal to each other.

\bigskip

Next, let us give the cases covering binary partitions.

\begin{corollary} \label{cf.23}
If $b_n(\mathfrak{2},k)$, is the number of unrestricted binary partitions of $k$ into non-negative powers of two no greater than $2^{n}$, then
    \begin{equation} \label{cf.24}
  \frac{1}{(1-q^1)(1-q^2)(1-q^4)...(1-q^{2^{n}})} = \sum_{k=0}^{\infty} b_n(\mathfrak{2},k) q^k
\end{equation}
\begin{equation*}
   = \cfrac{1}{1- \cfrac{q^1}{1 -\cfrac{q^2}{q^1 +q^2 - 2 q^1 q^2 -\cfrac{q^1 q^4}{q^2 + q^4 - 2 q^2 q^4 -\cfrac{\ddots}{\ddots
   \cfrac{q^{2^{n-2}} q^{2^{n}}}{q^{2^{n-1}} + q^{2^{n}} - 2q^{2^{n-1}}q^{2^{n}}}}}}}} \; .
\end{equation*}
\end{corollary}

The following distinct binary partitions example is completely solvable.

\begin{corollary} \label{cf.25}
If $p_n(\mathfrak{2D},k)$, is the number of binary partitions of $k$ into distinct non-negative powers of two no greater than $2^{n}$, then
    \begin{equation} \label{cf.26}
  (1+q^1)(1+q^2)(1+q^4)...(1+q^{2^{n}})= \frac{1-q^{2^{n+1}}}{1-q} = \sum_{k=0}^{2^{n+1}-1} p_n(\mathfrak{2D},k) q^k
\end{equation}
\begin{equation*}
   = \cfrac{1}{1- \cfrac{q^1}{1 + q^1 -\cfrac{(1+q^1)q^2}{q^1 + q^2 + q^1 q^2 -\cfrac{(1+q^2)q^4}{q^2 + q^4 + q^2 q^4 -\cfrac{\ddots}{\ddots \cfrac{(1+q^{2^{n-1\emph{}}})q^{2^{n}}}{q^{2^{n-1}}+q^{2^{n}}+q^{2^{n-1}}q^{2^{n}}}}}}}} \; .
\end{equation*}
\end{corollary}

Note that from (\ref{cf.26}) we have directly that
\begin{equation*}
  p_n(\mathfrak{2D},k) =
  \left\{
    \begin{array}{ll}
      1, & \hbox{when $0 \leq k < 2^{n+1}$;} \\
      0, & \hbox{when $k \geq 2^{n+1}$.}
    \end{array}
  \right.
\end{equation*}

The following distinct ternary partitions example is easily stated.

\begin{corollary} \label{cf.26a}
If $p_n(\mathfrak{3D},k)$, is the number of ternary partitions of $k$ into distinct non-negative powers of three no greater than $3^{n}$, then
    \begin{equation} \label{cf.26b}
  (1+q^1)(1+q^3)(1+q^9)...(1+q^{3^{n}})= \sum_{k=0}^{3^{n}-1} p_n(\mathfrak{3D},k) q^k
\end{equation}
\begin{equation*}
   = \cfrac{1}{1- \cfrac{q^1}{1 + q^1 -\cfrac{(1+q^1)q^3}{q^1 + q^3 + q^1 q^3 -\cfrac{(1+q^3)q^9}{q^3 + q^9 + q^3 q^9 -\cfrac{\ddots}{\ddots \cfrac{(1+q^{3^{n-1\emph{}}})q^{3^{n}}}{q^{3^{n-1}}+q^{3^{n}}+q^{3^{n-1}}q^{3^{n}}}}}}}} \; .
\end{equation*}
\end{corollary}

Note that from (\ref{cf.26b}) we have directly that
\begin{equation*}
  p_n(\mathfrak{3D},k) =
  \left\{
    \begin{array}{ll}
      1, & \hbox{for $0 \leq k < 3^{n+1}$; $k$ is a sum of distinct powers of 3.} \\
      0, & \hbox{for $0 \leq k < 3^{n+1}$; $k$ not a sum of distinct powers of 3.} \\
      0, & \hbox{for $k \geq 3^{n+1}$.}
    \end{array}
  \right.
\end{equation*}

Clearly this topic of Euler Continued Fractions applied to partition generating functions is an interesting elementary study for students, and a possible tool for researchers. The above results are old, and have probably been well-worked over time.

\section{Rogers-Ramanujan Continued Fractions for partition functions}

The fraction given here was mentioned by Ramanujan in his second letter to Hardy (see Adiga et al. \cite[p. xxviii]{cA1985}); namely

\begin{equation}  \label{cf.27}
  R(a,b)= 1+\cfrac{bq}{1+aq+ \cfrac{bq^2}{1+aq^2+\cfrac{bq^3}{1+aq^3+\cfrac{bq^4}{\ddots}}}} \; .
\end{equation}

However, these now famous continued fractions, as with the Rogers-Ramanujan identities, were first discovered in 1894 by Rogers (see \cite{lR1894}).
We define the functions $G(q)$ and $H(q)$ in the context of the Rogers–Ramanujan identities,
\begin{equation}  \label{cf.27a}
  G(q)= \sum_{n=0}^{\infty} \frac{q^{n^2}}{(1-q)(1-q^2)\cdots(1-q^n)} = \sum_{n=0}^{\infty} \frac{q^{n^2}}{(q:q)_n}
\end{equation}
\begin{equation}  \nonumber
  = \frac{1}{(q;q^5)(q^4;q^5)} = \prod_{n=1}^{\infty} \frac{1}{(1-q^{5n-4})(1-q^{5n-1})},
\end{equation}
and
\begin{equation}  \label{cf.27b}
  H(q)= \sum_{n=0}^{\infty} \frac{q^{n^2+n}}{(1-q)(1-q^2)\cdots(1-q^n)} = \sum_{n=0}^{\infty} \frac{q^{n^2+n}}{(q:q)_n}
\end{equation}
\begin{equation}  \nonumber
  = \frac{1}{(q^2;q^5)(q^3;q^5)} = \prod_{n=1}^{\infty} \frac{1}{(1-q^{5n-3})(1-q^{5n-2})}.
\end{equation}
The Rogers–Ramanujan continued fraction is then,
\begin{equation}  \label{cf.27c}
  R(q)= \frac{q^{\frac{11}{60}H(q)}}{q^{\frac{-1}{60}G(q)}}
      = q^{\frac{1}{5}} \prod_{n=1}^{\infty} \frac{(1-q^{5n-4})(1-q^{5n-1})}{(1-q^{5n-3})(1-q^{5n-2})}
\end{equation}
\begin{equation}  \nonumber
  = 1+\cfrac{q^{\frac{1}{5}}}{1+\cfrac{q}{1+\cfrac{q^2}{1+\cfrac{q^3}{\ddots}}}} \; .
\end{equation}

So, we note that $R(0, 1)$ leads us to the celebrated Rogers-Ramanujan continued fraction, which has been researched by many (see Andrews \cite[Chapter 7]{gA1976}, for example). In the course of analyzing identities from Ramanujan's Lost Notebook \cite{gA2018}, Andrews and Berndt have discussed the fraction $R(a, b)$, but mainly from the viewpoint of transformation formulas. Our emphasis here is on using (\ref{cf.27}) in a generalized approach to several partition identities, but there is a whole adjacent theory on particular values of these continued fractions determined from applying the theory of modular forms. Hence the examples, using $\varphi$ as the golden ratio $(\sqrt{5} +1)/2$,

\begin{equation}  \label{cf.28}
  \cfrac{e^{-\frac{-\pi}{5}}}{1+\cfrac{e^{-\pi}}{1+\cfrac{e^{-2\pi}}{1+\cfrac{e^{-3\pi}}{\ddots}}}}
   = \frac{1}{2}\varphi (\sqrt{5}-\varphi^{3/2}) (\sqrt[4]{5} + \varphi^{3/2}),
\end{equation}
\begin{equation}  \label{cf.29}
  \cfrac{e^{-\frac{-2\pi}{5}}}{1+\cfrac{e^{-2\pi}}{1+\cfrac{e^{-4\pi}}{1+\cfrac{e^{-6\pi}}{\ddots}}}}
   = \sqrt[4]{5} \varphi^{1/2} - \varphi,
\end{equation}
\begin{equation}  \label{cf.30}
  \cfrac{e^{-\frac{-4\pi}{5}}}{1+\cfrac{e^{-4\pi}}{1+\cfrac{e^{-8\pi}}{1+\cfrac{e^{-12\pi}}{\ddots}}}}
   = \frac{1}{2}\varphi (\sqrt{5}-\varphi^{3/2}) (-\sqrt[4]{5} + \varphi^{3/2}).
\end{equation}

So next we examine the continued fraction $R(a, b)$ of Ramanujan and consider various restricted partition functions. For further reading, a good reference is Alladi and Gordon \cite{kA1993}. We use the continued fraction to give results for several partition identities, some of which generalize results of Bressoud \cite{dB1979} and G\"{o}llnitz \cite{hG1967}. We also give a combinatorial interpretation for the coefficients in the power series expansion of the reciprocal $\frac{1}{R(-a, -b)}$, extending a result of Odlyzko and Wilf \cite{aO1988}. The full description of this approach would add several more pages to our work, but \cite{kA1993} covers all of this very nicely.

It turns out that \textbf{Lebesgue's identity} plays a major role in our analysis with respect to the numerators and denominators of the finite continued fractions we consider.
\begin{equation}  \label{cf.30a}
  \sum_{k \geq 0} \frac{q^{k(k+1)/2} \prod_{j=1}^{k} (1+bq^j)}{(1-q)(1-q^2)...(1-q^k)}
   = \prod_{m \geq 1}(1+bq^{2m})(1+q^m).
\end{equation}

It is known that Lebesgue's identity implies Ramanujan's fraction $R(a, b)$ has a product representation when $a = 1$. More precisely (\ref{cf.35}) and (\ref{cf.36}) (see below) yield
\begin{equation}  \label{cf.30b}
  1+\cfrac{bq}{1+q+ \cfrac{bq^2}{1+q^2+\cfrac{bq^3}{1+q^3+\cfrac{bq^4}{\ddots}}}} = \prod_{m=1}^{\infty} \frac{(1+bq^{2m-1})}{(1+bq^{2m})}.
\end{equation}
A neat case of (\ref{cf.30b}) is obtained from $q \mapsto q^2$ and $b \mapsto bq^{-1}$ so then
\begin{equation}  \label{cf.30c}
  1+\cfrac{bq}{1+q^2+ \cfrac{bq^3}{1+q^4+\cfrac{bq^5}{1+q^6+\cfrac{bq^7}{\ddots}}}} = \prod_{m=1}^{\infty} \frac{(1+bq^{4m-3})}{(1+bq^{4m-1})}.
\end{equation}

For a continued fraction $F$, let $P_n/Q_n$ denote its $n$th convergent, and suppose that $\lim_{n \rightarrow \infty} P_n = P$, $\lim_{n \rightarrow \infty} Q_n = Q$ in a suitable topology. We then say that $F$ has numerator $P$ and denominator $Q$, and write $P = F^N$, $Q=F^D$.
Consider the fraction

\begin{equation*}
  F(a,c)= 1+a+\cfrac{acq}{1+aq+ \cfrac{acq^2}{1+aq^2+\cfrac{acq^3}{1+aq^3+\cfrac{acq^4}{\ddots}}}} \; .
\end{equation*}

This can be written in the form

\begin{equation*}
  F(a,c)= \frac{f(a,c)}{f(aq,c)}, \quad \textmd{where} \quad f(a,c)=\sum_{k \geq 0} A_k q^k .
\end{equation*}

We now compute the coefficients $A_k=A_k(c, q)$, observing that $f(a, c)$ satisfies the recurrence
\begin{equation*}
  f(a,c)= (1 + a) f(aq,c)+ acq \, f(aq^2,c).
\end{equation*}

Therefore the coefficients $A_k$ satisfy
\begin{equation*}
  A_k = q^k \, A_k + q^{k-1} A_{k-1} \, q- cq^{2k-1}\, A_{k-1},
\end{equation*}
which is the same as
\begin{equation*}
  A_k = \frac{q^{k-1}(1+cq^k)}{(1-q^k)} A_{k-1}.
\end{equation*}

By iteration this yields
\begin{equation*}
  F(a,c)= \sum_{k \geq 0} \frac{a^k q^{\frac{k(k-1)}{2}}(-cq)_k}{(q)_k} .
\end{equation*}

Let $c = a^{-1}b$. Then
\begin{equation*}
  R(a,b)= \frac{f(a,a^{-1}b)}{f(aq,a^{-1}b)} - a
\end{equation*}
is Ramanujan's fraction (\ref{cf.27}).

\begin{lemma} \label{cf.31}
For the fraction $R(a, b)$, the numerator is
  \begin{equation} \label{cf.32}
  R^N(a,b)=  \sum_{k \geq 0} \frac{a^k q^{k(k+1)/2}(-a^{-1}b)_k}{(q)_k},
\end{equation}
and the denominator is
\begin{equation} \label{cf.33}
  R^D(a,b)=  \sum_{k \geq 0} \frac{a^k q^{k(k+1)/2}(-a^{-1}bq)_k}{(q)_k}.
\end{equation}
\end{lemma}

\textit{Proof}: The expansion (\ref{cf.33}) is an immediate consequence of
\begin{equation} \label{cf.34}
  R^D(a,b) =  f(aq, a^{-1}b).
\end{equation}

The expansion (\ref{cf.32}) is more complicated. To obtain it, observe that
\begin{eqnarray*}
  R^N(a,b) &=& f(a, a^{-1}b) - a \, f(aq, a^{-1}b) \\
           &=& \sum_{k \geq 0} \frac{a^k q^{k(k-1)/2}(-a^{-1}bq)_k}{(q)_k} -  \sum_{k \geq 0} \frac{a^{k+1} q^{k(k+1)/2}(-a^{-1}bq)_k}{(q)_k} \\
           &=& 1 +  \sum_{k \geq 0} \frac{a^{k+1} q^{k(k+1)/2}(-a^{-1}bq)_k}{(q)_k} \left( \frac{1+a{-1}bq^{k+1}}{1-q^{k+1}} -1 \right) \\
           &=& 1 +  \sum_{k \geq 0} \frac{a^{k+1} q^{(k+1)(k+2)/2}(-a^{-1}bq)_k (1- a^{-1}b)}{(q)_{k+1}}  \\
           &=& \sum_{k \geq 0} \frac{a^k q^{k(k+1)/2}(-a^{-1}b)_k}{(q)_k}
\end{eqnarray*}
as required. $ \quad \blacksquare$

Andrews (see \cite{gA1979a} and \cite{gA1981a}) considered the expansions in lemma \ref{cf.31} while discussing a transformation formula of Ramanujan \cite{sR1988} for $R(a, b)$. Our emphasis here is on the partition theorems that can be derived using $R(a, b)$, and for this the following lemma is crucial.

\begin{lemma}
  For the fraction $R(a, b)$, we also have the expansions
  \begin{equation} \label{cf.35}
  R^N(a,b)=  \sum_{i,j \geq 0} \frac{a^i b^j q^{(i^2+i)/2 + ij + j^2}}{(q)_i (q)_j},
\end{equation}
and the denominator is
\begin{equation} \label{cf.36}
  R^D(a,b)=  \sum_{i,j \geq 0} \frac{a^i b^j q^{(i^2+i)/2 + ij + j^2 + j}}{(q)_i (q)_j}.
\end{equation}
\end{lemma}

\textit{Proof}: To obtain (\ref{cf.35}) and (\ref{cf.36}) from (\ref{cf.33}) and (\ref{cf.34}) we use the $q$-binomial theorem,
\begin{equation*}
  (-z)_k =  \sum_{j=0}^{k} z^j q^{j(j-1)/2} \binom{k}{j}_q
\end{equation*}
with $z=a^{-1}b$ and $z=a^{-1}bq$. (See Campbell \cite{gC2019} for the $n$-space $q$-binomial theorem.) Therefore
\begin{eqnarray}
\nonumber   R^N(a,b) &=& \sum_{k \geq 0} \frac{a^k q^{k(k+1)/2}}{(q)_k}    \sum_{j=0}^{k} \frac{a^{-j} b^j q^{j(j-1)/2} (q)_k}{(q)_j (q)_{j-k}} \\
\nonumber            &=& \sum_{i,j \geq 0} \frac{a^{i} b^j q^{(i+j)(i+j+1)/2}}{(q)_i (q)_j},
\end{eqnarray}
where $i=k-j$; this is equivalent to (\ref{cf.33}). To obtain (\ref{cf.34}), observe that
\begin{equation} \label{cf.36a}
  R^D(a, b) = R^N(a, bq)
\end{equation}
by comparing (\ref{cf.35}) and (\ref{cf.36}).
\bigskip

The following two theorems relate successively to the numerator and the denominator of the fraction (\ref{cf.27}), so then to (\ref{cf.35}) and (\ref{cf.36}). For a proof of these see Alladi and Gordon \cite{kA1993}.

\begin{theorem}  \label{cf.37} \textbf{(Numerator)}

  Let $A^N(n; i,j)$ be the number of partitions of $n$ into $i+j$ distinct red parts and $j$ distinct blue parts such that one of the blue parts may
be zero and every blue part is $\leq i +j-1$.

Let $B^N(n; i,j)$ be the number of partitions of $n$ into $i$ distinct red parts and $j$ distinct non-consecutive blue parts such that every red part is $>j$.

Let $C^N(n; i,j)$ be the number of partitions of $n$ into $i$ red parts and $j$ blue parts such that all parts are distinct and after each blue part there is a gap of at least 2. Then
\begin{equation*}
  A^N(n; i,j) = B^N(n; i,j) = C^N(n; i,j).
\end{equation*}
\end{theorem}

\begin{theorem} \textbf{(Denominator)}

  Let $A^D(n; i,j)$ be as in $A^N(n; i,j)$ except that every blue part is $>0$ and $\leq i + j$.

  Let $B^D(n; i,j)$ be as in $B^N(n; i,j)$ except that part 1 cannot be blue.

  Let $C^D(n; i,j)$ be as in $C^N(n; i,j)$ except that part 1 cannot be blue. Then
\begin{equation*}
  A^D(n; i,j) = B^D(n; i,j) = C^D(n; i,j).
\end{equation*}
\end{theorem}

So reprising (\ref{cf.30c}) namely
\begin{equation*}
  1+\cfrac{bq}{1+q^2+ \cfrac{bq^3}{1+q^4+\cfrac{bq^5}{1+q^6+\cfrac{bq^7}{\ddots}}}} = \prod_{m=1}^{\infty} \frac{(1+bq^{4m-3})}{(1+bq^{4m-1})}),
\end{equation*}

we have interesting cancellations in numerator-denominator equations. That is, the numerator is given by
\begin{eqnarray*}
  \sum_{k \geq 0} \frac{q^{k(k+1)}(-bq^{-1};q^2)_k}{(q^2;q^2)_k} &=& \prod_{m=1}^{\infty} (1+bq^{4m-3})(1+q^{2m}) \\
                  &=& \prod_{m=1}^{\infty} (1+bq^{4m-3})(1+q^{4m-2})(1+q^{4m})
\end{eqnarray*}
and the denominator is given by
\begin{equation*}
  \sum_{k \geq 0} \frac{q^{k(k+1)}(-bq;q^2)_k}{(q^2;q^2)_k} = \prod_{m=1}^{\infty} (1+bq^{4m-1})(1+q^{4m-2})(1+q^{4m})
\end{equation*}
with right sides having common factors that eliminate.

This leads in particular to the continued fraction identity
\begin{equation}  \label{cf.37a}
  1+\cfrac{q}{1+q^2+ \cfrac{q^3}{1+q^4+\cfrac{q^5}{1+q^6+\cfrac{q^7}{\ddots}}}}
   = \frac{\prod_{j \equiv 2,3,7 \, (\textmd{mod} 8)}(1-q^j)}{\prod_{j \equiv 1,5,6 \, (\textmd{mod} 8)}(1-q^j)}.
\end{equation}
G\"{o}11nitz \cite{hG1967} states similar results, but (\ref{cf.37}) seems to have escaped attention.
There is a continued fraction identity due to Gordon \cite{bG1965} and G\"{o}11nitz \cite{hG1967} which looks very similar to (\ref{cf.37a}), namely
\begin{equation}  \label{cf.38}
  1+q+\cfrac{q^2}{1+q^3+ \cfrac{q^4}{1+q^5+\cfrac{q^4}{1+q^7+\cfrac{bq^6}{\ddots}}}}
   = \frac{\prod_{j \equiv 3,4,5 \, (\textmd{mod} 8)}(1-q^j)}{\prod_{j \equiv 1,4,7 \, (\textmd{mod} 8)}(1-q^j)}.
\end{equation}
However, this result first appears in Alladi and Gordon \cite{kA1993} almost 30 years after (\ref{cf.37}).

\section{Ramanujan's three parameter continued fraction}

Ramanujan \cite{sR1957} obtained in addition to (\ref{cf.27}), the following continued fraction with three parameters $a$, $b$, $q$ which has also a product representation
\begin{equation}  \label{cf.39}
  1-ab+\cfrac{(a-bq)(b-aq)}{(1-ab)(1+q^2)+ \cfrac{(a-bq^3)(b-aq^3)}{(1-ab)(1+q^4)+\cfrac{(a-bq^5)(b-aq^5)}{(1-ab)(1+q^6)+\cfrac{(a-bq^7)(b-aq^7)}{\ddots}}}}
\end{equation}
\begin{equation*}
   = \prod_{m=1}^{\infty} \frac{(1+a^2 q^{4m-3})(1+b^2 q^{4m-3})}{(1+a^2 q^{4m-1})(1+b^2 q^{4m-1})}.
\end{equation*}
This was proved only in 1985 by the reviewers of Chapter 16 of Ramanujan's Second Notebook \cite{cA1985}, 65 years after Ramanujan's death. If we put $a=0$ and replace $b^2$ by $-b$ in (\ref{cf.39}), we get (\ref{cf.30c}). It seems there is still scope to study the combinatorial properties of the coefficients in the power series expansion of this fraction.

\end{document}